\documentclass{statsoc}
\usepackage{amsmath}
\usepackage[dvips]{graphicx}

\title[Negative Surveys]{Negative Surveys}
\author[Fernando Esponda]{Fernando Esponda}
\address{Yale University,
   Computer Science Department,
  New Haven CT 06520, USA}
\email{fernando.esponda@yale.edu}
\begin{document}

\maketitle
\begin{abstract}

In this paper we propose a strategy for administering a survey that is
mindful of sensitive data and individual privacy. The survey in
question seeks to estimate the population proportions of a sensitive,
polychotomous variable and does not depend on anonymity, cryptography,
or in legal guarantees for its privacy preserving properties.

Our technique, called {\em Negative Surveys}, presents interviewees
with a question and $t$ possible answers, and asks participants to
eliminate one of $t-1$ alternatives at random.  The method is closely
related to randomized response techniques (RRTs) in that both rely on
a random component to preserve privacy; however, while RRTs require
respondents to choose among questions and give an answer, negative
surveys ask them to choose between possible answers to a single
question. This distinction has important consequences for the privacy,
methodology, and reliability of our scheme.

In the course of the paper we quantify the amount of information
surrendered by an interviewee, elaborate on how to estimate the
desired population proportions, and discuss the properties of our
method at length. We also introduce a specific setup that requires a
single coin as a randomizing device, and that limits the amount of
information each respondent is exposed to by presenting to them only a
subset of the question's alternatives.


\keywords{Privacy preserving surveys; Estimating population proportions;
  Randomized response techniques; Sensitive attribute; Polychotomous variable}

\end{abstract}
\section{Introduction}\label{intro}


Surveys are an indispensable tool for learning the characteristics of
a population, from matters of public health, such as drug use
frequency, to gauging public opinion on issues like abortion and
homosexuality, to electing the leaders of a democratic country.  Their
reliability is dependent on having ample participation and an unbiased
sample. On occasion, however, they require disclosing sensitive or
otherwise private data compelling some interviewees to give faulty
answers and discouraging others from participating altogether.  There
is seldom a big incentive to answer a survey, and when its questions
are potentially stigmatizing, special care should be taken to protect
respondent privacy and promote participation.



In this paper, we propose {\em Negative surveys} as a technique for
conducting surveys that is mindful of participant's privacy.  Negative
surveys allow participants to keep the target datum undisclosed by
asking them, instead, to make a series of decisions with the datum in
mind.  In this way, the frequency of drug use can be calculated
without respondents admitting to using any drugs, the popular opinion
on abortion can be measured without asking for anyone's specific
position, and an election can be run without any of the voters
explicitly stating their preference.  Our objective is that, by
providing transparent privacy guarantees, studies using our scheme
will have greater participation and more accurate responses, and will,
therefore, be more reliable.

In what follows, we review some of the work related to our own,
specify the kind of survey under study and the negative survey
technique, and discuss the privacy of our method by analyzing the
amount of information gained by each questionnaire. We then look at
how a negative survey is applied, how the sought data is extracted
from it, and, finally, conclude by summarizing its characteristics.


\subsection{Related Work}
\cite{mythesis} studies the idea of depicting a data set
$DB$ by alternatively storing every datum {\em not} in $DB$---the
universe of possible data items is assumed to be finite. This is
accomplished by introducing a data compaction scheme and a series of
algorithms that allow the complement set to be created and stored
efficiently.  One property of this construction, called {\em negative
database}, is its potential for restricting the kinds of inferences
that can be drawn from the data: given a negative database, it is a
difficult problem to recover the original dataset $DB$; yet, answers
for a certain, limited type of queries can be obtained efficiently.

Other properties or representing data negatively are outlined in
\cite{mythesis}; including the suggestion that this viewpoint may be
useful, not only to protect stored data, but also, as a paradigm to
enhance the privacy of data collection.  It is this proposition that
inspired the current work.

A suit of techniques that share the same motivation as our proposal---
to protect privacy, promote participation, and increase survey
accuracy---and which have a similar procedure for conducting surveys
are known as {\em Randomized response techniques}.

Randomized response techniques (RRT) were introduced in
\cite{warner65}. The original model sets out to estimate the
proportion of a population that belongs to a particular, stigmatizing
group $A$. It offers participants two possible questions:
\begin{itemize}
\item[$Q1$]: Do you belong to group $A$?
\item[$Q2$]: Do you belong to group $B$?
\end{itemize}
where $A$ and $B$ are exhaustive and mutually exclusive groups, i.e.,
$B=\bar{A}$; for example, if $A$ represents the group of people that
have had sex with a minor; $B$ would stand for the people that
haven't. Only one of the two questions is to be answered.  The
question must be selected privately using a randomizing device
provided by the interviewer and should remain undisclosed at all
times. The only information surrendered by the participant is a Yes or
No answer, not the question being answered, not the outcome of the
randomizing device.  In this way, the interviewee avoids disclosing
which group he belongs to, yet provides sufficient information to
estimate the desired proportions---the Yes and No answers of all the
respondents in the sample along with the known characteristics of the
randomizing device, are enough to estimate the proportion of
population members in each group (see
\cite{warner65,fox86,chaud88,mangat94,gjestvang05}).

A further refinement of this approach relaxes the need for groups $A$
and $B$ to be exhaustive and mutually exclusive and substitutes the
second question for something less sensitive. For example, $Q1$ could
read ``Have you had sex with a minor?'' and $Q2$: ``Do you belong to
the YMCA?''  This variant is know as the unrelated question or
paired-alternative method \cite{horvitz67,moors71,chaud88}.

The RRT methods discussed so far are designed for dichotomous
populations. Negative surveys, on the other hand, are relevant only
when there are more than two categories in which the population can be
divided. \cite{Abul67} generalizes Warner's model to polychotomous
populations by employing several independent
samples. \cite{bourke73,bourke76} propose a different scheme that
necessitates only a single sample: categories are numbered 1 through
$t$, participants disclose which category they belong to with
probability $p$ or, with probability $1-p$, choose a number between 1
and $t$; each number is selected with probability $p_1 \ldots p_t$,
where $\sum p_i=1-p$ (see \cite{chaud88} for more detail and
\cite{kimwarde05} for a more recent example).

Finally, it is worth pointing out the existence of mechanisms for
conducting direct response surveys privately. For instance, anonymity
schemes can be used to conceal the identity and input of individual
participants, e.g.,  \cite{sudman74}'s self administered questionnaires,
Web and e-mail anonymous surveys, and cryptographically based surveys
as in \cite{feigenbaum04}. Also, legal guarantees can be set that
safeguard respondent privacy.  However, these methods may not always
lead to the desired participation level since respondents still need
to answer a sensitive question---the guarantees offered by the study
need to be understood to be trusted, and trust must be put on higher
authorities not to circumvent the promise. Further, some of these
techniques require setups that are not always available for a study,
such as the use of computers, and some, like anonymous surveys, have
additional shortcomings when it comes to verifying their results or
conducting longitudinal studies.  \cite{fox86} discuss the drawbacks of
some of these approaches at length.


\section{Negative Surveys}
The type of survey we consider consists of a questionnaire with a
single question (or statement) and $t$ categories
$\{X_1,X_2,\ldots,X_t\}$ from which to choose an answer (or
alternative).  The survey is administered to a sample of $n$
individuals drawn uniformly at random with replacement from the
population.

We refer to it as a {\it Positive survey} or as a {\it Direct response
survey } when the subjects are asked to reveal which category they
belong to. We call it a {\it negative survey} when the requirement is
to disclose a category (a single one, for the current work) to which
they {\em do not} belong---a negative questionnaire can be obtained by
simply negating the question of a positive questionnaire.  The
categories in the direct response survey are exhaustive and mutually
exclusive---one and only one option is true; in the negative survey,
one and only one category is false for a particular individual.  The
object of both versions of the survey is to estimate the proportions
of the population that belong to each category.

Take, for example, a direct response, salary survey:\\
\begin{center}
\begin{minipage}[h]{0.7\linewidth}
I earn:
\begin{itemize}
\item[{\bf[   ]}] Less that 30,000 dollars a year
\item[{\bf[   ]}] Between 30,000 and 60,000 dollars a year
\item[{\bf[   ]}] More that 60,000 dollars a year
\end{itemize}
\end{minipage}
\end{center}

The negative version would read:
\begin{center}
\begin{minipage}[h]{0.7\linewidth}
I {\bf do not} earn:
\begin{itemize}
\item[{\bf[   ]}] Less that 30,000 dollars a year
\item[{\bf[   ]}] Between 30,000 and 60,000 dollars a year
\item[{\bf[   ]}] More that 60,000 dollars a year
\end{itemize}
\end{minipage}
\end{center}

If the positive version of the survey is being answered by an
individual whose income is 20,000 dollars, the first option must be
chosen. Alternatively, if the same person is answering the negative
version of the survey, one of the last two options should be selected.
Next we look more closely at the amount of information that is being
surrendered in both cases.

\subsection{Privacy Preserved}\label{information}


It is intuitively clear that the amount of information required by a
negative questionnaire is inferior to what is asked for in its direct
response version, at least when the query has more than two options.
We formalize this notion and show that, indeed, the information
required for a negative survey is at most that of its positive
counterpart.

Using Shannon's uncertainty measure, \cite{shannon48}, the amount of
information gained from a positive questionnaire in which categories
are exhaustive and mutually exclusive can be written as:
\begin{equation}\label{posinfo}
-\sum_i p_i \log p_i
\end{equation}
where $p_i$ is the probability that option $X_i$ is true and $t$ is
the number of categories in the questionnaire. The maximum amount is
obtained when all options are equally likely.

Now consider the information gained from applying a negative
questionnaire in which only one option, $X_s$, is selected by the
respondent. We compute this quantity as the difference in information
of two positive questionnaires: the information obtained by the
positive version of the questionnaire (given in Eq. \ref{posinfo}),
minus the information gained from the same questionnaire once $X_s$ is
no longer an option.

\begin{equation}\label{neginfo}
-\sum_i p_i \log p_i +\sum_{i\neq s} P(X_i=T|X_s=F) \log P(X_i=T|X_s=F)
\end{equation}
where $P(X_i=T|X_s=F)$ is the probability that category $i$ is true in
a direct response survey after $X_s$ has been removed as an option,
i.e., after finding out it is false.

It is easy to see from the above expressions that the information
gained from a negative questionnaire is at most the quantity obtained
from its positive counterpart.

%

\section{Estimating Proportions Using  Negative Input}\label{Estimating}
In the previous sections, we explained negative surveys and discussed
how their application increases the interviewees' privacy by requiring
less information than their positive counterpart.  Negative surveys
ask respondents to choose one of the $t-1$ options that truthfully
answer the question before them---all choices except one are true for
a specific individual when surveyed ``negatively''.  However, we are
after the proportions of the population that ``positively'' belong to
each of the $t$ categories: a particular interviewee positively
belongs to one and only one category.
In this section, we show how to estimate these values along with their
corresponding measures of variation.  The analysis follows a similar
reasoning---albeit different in the details---as the one used for
randomized response techniques, particularly as shown in
\cite{chaud88} for vector responses. We therefore adopt the same
notation.

Let $p_{i,j}$ ($1\leq i,j \leq t$) be the probability that option
$X_i$ is chosen given that a respondent positively belongs to $X_j$,
and $\sum_i p_{i,j}=1$.  Let $\pi_i$ denote the proportion of the
population that positively belongs to category $i$ with
$\sum_i\pi_i=1$. Then, the probability of selecting $X_i$ is given by:
\begin{equation} \label{prob-lambda}
\lambda_i=\sum_j p_{i,j} \pi_j 
\end{equation}
Let $P$ denote the matrix of $p_{i,j}$'s:
\begin{equation*}
P=\left [
\begin{array} {cccc}
 p_{1,1} &+ p_{1,2} &+ \cdots &+ p_{1,t}\\
 p_{2,1} &+ p_{2,2} &+ \cdots &+ p_{2,t}\\
\vdots   &  &   \ddots &\\
 p_{t,1} &+ p_{t,2} &+ \cdots &+ p_{t,t} 
\end{array} \right ]
\end{equation*}
 where $\sum_i p_{i,j}=1$ and $p_{i,i}=0$; let
 $\boldsymbol{\pi}=(\pi_1,\ldots,\pi_t)'$ and
 $\boldsymbol{\lambda}=(\lambda_1,\ldots,\lambda_t)'$. The probability
 of responses for each category is written in matrix notation as:

\begin{equation}
P \boldsymbol{\pi}=\boldsymbol{\lambda} 
\end{equation}

Let $n_i$ be the observed frequency of category $X_i$ ($1\leq i\leq
t$) obtained from the application of a negative survey to $n$
individuals.  Observe that $n_i$ is binomially distributed with
parameters $n$ and $\lambda_i$.  An unbiased estimator of $\lambda_i$
is $\hat{\lambda_i}=\frac{n_i}{n}$, and, provided $P$ is non-singular,
an unbiased estimator of $\boldsymbol{\pi}$ is given by:
\begin{equation}
\boldsymbol{\hat{\pi}}=P^{-1}\boldsymbol{\hat{\lambda}}
\end{equation}
where 
$\boldsymbol{\hat{\pi}}=(\hat{\pi_1},\ldots,\hat{\pi_t})'$ and
$\boldsymbol{\hat{\lambda}}=(\hat{\lambda_1},\ldots,\hat{\lambda_t})'$.

An unbiased estimator for the variance and covariance of
$\boldsymbol{\hat{\pi}}$ is computed as follows:
Let $\boldsymbol{\hat{\lambda}^d}$ be a diagonal matrix where entry
$(i,i)$ is equal to the $i^{th}$ element of
$\boldsymbol{\hat{\lambda}}$; the estimated covariance of
$\boldsymbol{\hat{\pi}}$ is written as:

\begin{equation}
\hat{cov}(\hat{\boldsymbol{\pi}})=\frac{1}{n-1}P^{-1}(\boldsymbol{\hat{\lambda}^d}-\boldsymbol{\hat{\lambda}}\boldsymbol{\hat{\lambda}'})P^{'-1}
\end{equation}

The accuracy of the resulting $\pi_i$'s is dependent on an adequate
sampling of the $n_i$'s, as would also be the case in a positive
survey, and on having a good estimate of the $p_{i,j}{'s}$.  Knowing
how individuals choose an option is the extra information needed to
estimate the desired proportions while keeping personal data
concealed. Insight may come from knowledge about the behavior of the
population---factors like name recognition may bias an individual to
select a category---from data gathered in previous surveys, or from
employing a specific design in the administration of the survey. This
last alternative is discussed further in the next section.


\subsection{How to Choose an Option}\label{HowtoChoose}

In this section we propose a scheme intended to reduce the impact of
unknown biases by automating, and hence, predetermining part of the
decision process used to elect an answer.

One way to determine how respondents select a category is by
instructing them on how to choose among available options. For
instance, a simple, straightforward design gives each category an
equal chance of being selected:
\begin{equation*}
P=\left [
\begin{array} {llll}
 0 &+ \frac{1}{t-1} &+ \cdots &+ \frac{1}{t-1}\\
 \frac{1}{t-1}&+ 0 &+ \cdots &+ \frac{1}{t-1}\\
\vdots   &  & \quad  \ddots &\\
 \frac{1}{t-1} &+ \frac{1}{t-1} &+ \cdots &+ 0 
\end{array} \right ]
\end{equation*}

In this scenario the probability of option $X_i$ being chosen is:
\begin{equation}
\lambda_i=\frac{1}{t-1}(1-\pi_i) 
\end{equation}
following \ref{prob-lambda} and noting that $\sum_i\pi_i=1$. An
unbiased estimator of $\pi_i$ is given by:
\begin{equation} \label{eq-samep-for-pi}
\hat{\pi_i}=1-(t-1)\hat{\lambda}_i 
\end{equation}
where $\hat{\lambda_i}=n_i/n$ is an unbiased estimator for
$\lambda_i$. Similarly, an u.e. for the variance of $\hat{\pi_i}$:
\begin{equation}\label{eq-samep-for-var}
\hat{var}(\hat{\pi_i})=\frac{(t-1)^2}{n-1}\hat{\lambda}_i (1-\hat{\lambda}_i)
\end{equation}
and the covariance of  $\hat{\pi_i}$ and  $\hat{\pi_j}$:
\begin{equation}\label{eq-samep-for-cov}
\hat{cov}(\hat{\pi_i},\hat{\pi_j})=-\frac{(t-1)^2}{n-1}\hat{\lambda}_i
\hat{\lambda}_j
\end{equation}

With this scheme respondents provided with a fair, $t-1$ sided, die
can select an answer by privately obtaining a value $m$, and choosing,
for instance, the $m^{th}$ true option from the top, skipping over the
false category if needed.  One difficulty of the approach is, in fact,
having a $t-1$ sided die readily available, as during a phone survey.
This is compounded when asking several questions, each with a
different number of categories.  In the following section we propose a
design that addresses this point and illustrates other important
properties of gathering information with a negative survey.


\subsubsection{Two-option Survey Scheme}

It was earlier discussed how it is essential to know how respondents
choose, on average, from the available options, and considered using a
$t-1$ sided die to determine their selection. However, this has the
inconvenience of necessitating a customized device for each question
(unless an general purpose random number generator is at hand).  Next,
we present a scheme that resolves this issue and reduces the impact of
unknown biases (arising from an incorrect compliance with the survey
instructions) by automating part of the decision process.

Direct response surveys as well as randomized response techniques
require presenting all categories to the interviewee---for only one
choice is true. Conversely, in a negative survey all options except
one are true; consequently, only a subset needs to be evaluated by the
respondent.  With this in mind, our scheme preselects, uniformly at
random, a subset of the categories (two or more) for each of the
individuals questioned.

Consider the case where each subject is presented with a question and
two options; if both options prove true, one should be selected at
random---the interviewee privately tosses a fair coin prior to reading
the question and picks the first category if heads, the second if
tails---otherwise, when only one option is true, the true option must
be selected. Note that the setup only makes sense when the original
survey doesn't include categories of the type ``None of the above''.

The probability of choosing $X_i$ in this setup, is the probability of
it being presented in the questionnaire times the probability that it
is selected by the interviewee.  We can analyze this as the sum of two
terms \footnote{The probability of $X_i$ being selected if it does not
answer the question truthfully is considered to be zero and therefore
omitted.}: when it is presented alongside the only false option
($X_j$), in which case it will be selected, and when it is paired with
another true alternative and a choice must be made between them:

\begin{equation}\label{probXi}
p_{i,j}=\frac{2}{t(t-1)}+\sum_k \frac{2}{t(t-1)}P(X_i|X_k=T), \qquad
(k\not =j)
\end{equation}
where $P(X_i|X_k=T)$ denotes the probability of $X_i$ being chosen
given that it is presented together with another true option $X_k$.

\noindent
We obtain a bound for selecting $X_i$, when it is a true alternative,
by setting $P(X_i|X_k=T)$ to zero and one respectively:
\begin{equation}
p_{i,j} \in \left [\frac{2}{t(t-1)},\frac{2(t-2)}{t(t-1)} \right ]
\end{equation}
Consider that if all categories are presented to each respondent and
respondents ignore the outcome of the randomizing device, there might
be a particular alternative that never gets chosen (true also if the
subset shown to the interviewee has more than two categories), or an
option that always does (except when it proves false to the
interviewee). Using the two-option survey scheme, this possible error
is greatly reduced.

Finally, if no additional information is available regarding how
subjects choose between true alternatives, we assume
$P(X_i|X_k=T)=\frac{1}{2}$ for all $i \not=k$, and, therefore, that
$p_{i,j}=\frac{1}{t-1}$ for all $i \not=j$.
The sought after proportions, $\pi_i's$, are computed according to
Eq. \ref{eq-samep-for-pi} and the variances and covariances according
to Eq. \ref{eq-samep-for-var} and \ref{eq-samep-for-cov}
respectively.

Three interesting features of this approach are worth noting:
first, the use of a single randomizing device for every question
independent of the number of categories (for the above scenario, the
use of a coin instead of a $t-1$ sided die); second, the ability to
conduct a survey without disclosing all of the options to individual
respondents; and third, that the error in estimating the $p_{i,j}{'s}$
is bounded even if the coin is used improperly.





\section{Conclusion}

Survey accuracy depends on minimizing the incidence of nonrespondents
and on the honest participation of respondents. In studies where
interviewees are required to answer sensitive questions, care must be
taken to avert these difficulties by ensuring respondent privacy.
In this paper we presented a method for administering a questionnaire
that safeguards interviewee's privacy.  The survey in question seeks
to estimate the population frequencies of a polychotomous variable; it
consists of a single (potentially sensitive) question and $t$ options
from which to choose an answer.  Its privacy preserving properties do
not rely on anonymity, cryptography or on any legal contracts, but
rather on participants not revealing their true answer to the
survey's query---respondents are only required to discard, with a
known probability distribution, some of the categories that do {\em
not} answer the question for them.  This information is enough to
estimate the population proportions of the variable under study; yet,
insufficient to ascribe a sensitive datum to a particular individual.
We call the method Negative Surveys.

Negative surveys are closely related to randomized response techniques
(RRTs): both aim at conducting private surveys and both rely on the
(secret) use of a randomizing device to answer questionnaires.
One key distinction, however, is that in RRTs participants use the
device to choose among questions, at least one of which is sensitive;
while in negative surveys, they use it to choose among answers,
avoiding the problem of selecting the proper alternative question
altogether.  Also, with RRTs some subjects will be selected, by the
randomizing device, to answer the potentially stigmatizing
question. It might still be problematic for them to participate, as
the question remains sensitive and answering it demands a measure of
trust on the surveying scheme and on the surveyors---\cite{fox86} cite
a study in which the randomizing device was rigged in order to study
respondent behavior; a similar practice could be used to subvert their
privacy.  Negative surveys never prompt respondents to answer a
sensitive query directly.

We also presented a special setup for negative surveys that reduces
the complexity of the randomizing device to a simple, fair coin
revealing an important characteristic of our method: an interviewee
does not need to contemplate all of the question's potential answers
to pick his own, furnishing a level of secrecy to the survey itself
and providing robustness against the non-observance of questionnaire
instructions (\cite{ambainis98} discuss cryptographic techniques to
avoid cheating in RRT).

We expect that the privacy of a negative survey, its comprehensible
guarantees, and robustness will increase the level of cooperation and
accuracy in topic-sensitive studies.





\bibliographystyle{chicago}
\bibliography{$HOME/bibs/survey}

\end{document}